\documentclass[oneside]{article}
\usepackage{amsmath, enumerate, graphicx, fancyheadings, multirow, amssymb, array}
\usepackage[textwidth=5.5in, textheight=8in, footskip=0.45in, voffset=0.4in]{geometry}
\allowdisplaybreaks

\begin{document}
%
%
\pagestyle{fancy}
\lhead{TechS Vidya e-Journal of Research \\ ISSN 2322-0791 \\ 
Vol. 1, (2012-2013), pp.84-87.}
\rhead{{\small Distributions of the same type: non-equivalence of definitions $\ldots$}}
%
%
\ 

\begin{center}
\Large\bfseries 
Distributions of the same type: \\ 
non-equivalence of definitions in the discrete case 
\end{center}
%
%
\ 

\begin{center}
\large\bfseries S. Satheesh
\end{center}
%
%
\begin{center}
Department of Applied Sciences  \\ 
Vidya Academy of Science \& Technology, Thrissur, India \\
e-mail:ssatheesh1963@yahoo.co.in
\end{center}
%
%
\begin{center}
{\bf Abstract}
\end{center}
\begin{quote}
Distributions of the same type can be discussed in terms of distribution functions as well as their integral transforms. For continuous distributions they are equivalent. In this note it is shown that it is not so in the discrete case.
\begin{center}

\end{center}
\end{quote}
%
%
\label{SatheeshStart}
\section{Introduction}
Distributions of the same type have roles when we need to consider change of scale of a distribution for better fit to the given data (or in general scale families of distributions), standardization of distributions and in the area of stability problems in stochastic models. Distributions of the same type can be discussed in terms of distribution functions (d.f.) as well as their integral transforms (characteristic functions (CF), Laplace transforms (LT) or probability generating functions (PGF)). When we study summation stability of distributions it is convenient to consider distributions of the same type in terms of their integral transforms depending on their  support, that is,  CF's for distributions on $\mathbb R =(-\infty, \infty)$, LT's on for distributions on $\mathbb R_+=(0,\infty)$, 
and PGFs for those on $\mathbb Z_+ = \{0,1,2,\ldots\}$. When we study stability of extremes of distributions we will consider distributions of the same type in terms of their d.f.'s or survival functions. 

Let $X$   be a continuous random variable (r.v.) with d.f. $F(x)$, $x\in \mathbb R$ and CF $\phi(t)$, $t\in \mathbb R$. A change of scale to  $X$  is represented by $cX$, $c>0$,   and the corresponding d.f. and CF are given by 
\begin{align*}
G(x) & = P(cX\le x) \\
& = P\left(X\le \frac{x}{c}\right)\\
& = F\left( \frac{x}{c}\right)\\
\psi(t) & = E\left(e^{ictX}\right) \\
& = \phi(ct)
\end{align*}
Thus there is one-to-one correspondence between $F\left( \frac{x}{c}\right)$   and $\phi(ct)$. 
\section{Distributions of the same type: discrete case}
Now consider distributions of the same type for those on  $\mathbb Z_+$. The first definition is derived from the LT of the continuous analogue. 
\paragraph{Definition 2.1} (see [1])
{\em Two PGF's $Q_1(s)$  and $Q_2(s)$   are of the same type if  for all $s\in(0,1)$ and some $\alpha \in (0,1)$ we have
$$
Q_1(s)=Q_2(1-\alpha +\alpha s).
$$
}

The d.f. of a mixture of geometric distributions on $\mathbb Z_+$ has the general form
$$
F(k) = P(X<k) = 1- m(k), \quad k=0,1,\ldots
$$
where $\{m(k)\}$ is the moment sequence of the mixing distribution.
Further $\{m(k)\}$ is also the sequence of realizations of a LT $m(s)$, $s >0$, at the non-negative integral values of $s$.  Since $ m(\alpha s)$ also is a LT for a constant $\alpha >0$, we can define another d.f.  by
$$
G(k)=1-m(\alpha k),\quad k=0,1,\ldots
$$
Writing 
$$
F_\alpha(k)=1-m(\alpha k)
$$
we have
\begin{equation}
G(k) = F_\alpha(k)\quad\text{for all}\quad k=0,1,\ldots
\end{equation}
Thus we have a scale family $\{ F_\alpha(k): \alpha>0\}$ of d.f.'s on $\mathbb Z_+$ and $F(k)$ is a member of it. The existence of such a family of d.f.'s justifies the following definition.

\paragraph{Definition 2.2} (see [2])
{\em Two d.f.'s $F(k)$ and $G(k)$ on $\mathbb Z_+$ are of the same type if equation (1) is satisfied for some $\alpha>0$. }

\paragraph{Example 2.1}

Let $X$ has a geometric ($p$) distribution on $\mathbb Z_+$. Then its d.f.  is 
$$
F(k)=1-q^k,\quad k\in \mathbb Z_+
$$
and PGF is
$$
Q_X(s) = \frac{p}{1-qs}, \quad q=1-p.
$$
Now following Definition 2.2 for some $\alpha\in (0,1)$ (so that Definition 2.1 is also accommodated) consider the r.v.  $ Y$ with d.f. 
$$
G(k) = 1-q^{\alpha k},\quad k=0,1,2,\ldots
$$
Setting
$$
q=\frac{1}{4}, \quad \alpha = \frac{1}{2}
$$
we have
$$
q^\alpha = \sqrt{\frac{1}{4}}=\frac{1}{2}.
$$
Further we have
\begin{align*}
Q_X(s) & = \frac{3}{4-s} \\
Q_Y(s) & = \frac{1}{2s} \\
Q_X\left(1-\frac{1}{2}+\frac{s}{2} \right) & = \frac{6}{7-s} \\
Q_Y\left( \frac{1}{2} + \frac{s}{2} \right) & = \frac{2}{3-s}
\end{align*}
Thus we have
\begin{align*}
Q_X(s) & \ne Q_Y\left(\frac{1}{2}+\frac{s}{2}\right), \\
Q_Y(s) & \ne Q_X\left(\frac{1}{2}+\frac{s}{2}\right),
\end{align*}
considering both the possibilities. 

Hence the two definitions are not equivalent. This example has appeared in [3].
\end{document}